\documentclass[11pt]{article}
\usepackage{epic,latexsym,amssymb}
\usepackage{amsfonts}
\usepackage{amscd}
\usepackage{amsmath}
\usepackage{graphicx}
\usepackage{color}
\usepackage{caption,
caption}
\usepackage{tikz}

\usepackage{lineno}

\usetikzlibrary{arrows.meta}
\usepackage{float}

\newtheorem{thm}{Theorem}

\newtheorem{prop}[thm]{Proposition}
\newtheorem{lem}[thm]{Lemma}

\newtheorem{cor}[thm]{Corollary}
\newtheorem{prob}{Problem}

\newcommand{\proof}{\noindent\textbf{Proof. }}

\newcommand{\qed}{$\Box$}

\newcommand{\QEDmark}{\mbox{\textsc{qed}}}
\newcommand{\proofStarter}[1]{\textsc{#1} }

\def\vertex(#1){\put(#1){\circle*{2}}}
\def\vertexo(#1){\put(#1){\circle{2}}}
\def\vert(#1){\put(#1){\circle*{1.5}}}
\def\verto(#1){\put(#1){\circle{1.5}}}
\def\lab(#1)#2{\put(#1){\makebox(0,0)[c]{#2}}}
\newcommand{\mof}{\mathrm{mof}}
\newcommand{\MOF}{\mathrm{MOF}}

\definecolor{DarkGreen}{rgb}{0.2, 0.6, 0.3}

\definecolor{electricindigo}{rgb}{0.44, 0.0, 1.0}

\let\oldenumerate\enumerate
\renewcommand{\enumerate}{
  \oldenumerate
  \setlength{\itemsep}{0.5pt}
  \setlength{\parskip}{0pt}
  \setlength{\parsep}{0pt}
}

\title{Maximum oriented forcing number for complete graphs}
\author{Yair Caro\thanks{Dept. of Math and Physics, University of Haifa-Oranim, Israel} and Ryan Pepper\thanks{Dept. of Math and Statistics, University of Houston--Downtown, USA -- corresponding author}}

\begin{document}
\maketitle
\begin{abstract}
The \emph{maximum oriented $k$-forcing number} of a simple graph $G$, written $\MOF_k(G)$, is the maximum \emph{directed $k$-forcing number} among all orientations of $G$. This invariant was recently introduced by Caro, Davila and Pepper in~\cite{CaroDavilaPepper}, and in the current paper we study the special case where $G$ is the complete graph with order $n$, denoted $K_n$. While $\MOF_k(G)$ is an invariant for the underlying simple graph $G$, $\MOF_k(K_n)$ can also be interpreted as an interesting property for tournaments. Our main results further focus on the case when $k=1$. These include a lower bound on $\MOF(K_n)$ of roughly $\frac{3}{4}n$, and for $n\ge 2$, a lower bound of $n - \frac{2n}{\log_2(n)}$. Along the way, we also consider various lower bounds on the maximum oriented $k$-forcing number for the closely related complete $q$-partite graphs.
\end{abstract}
{\small \textbf{Keywords:} tournaments, maximum oriented forcing, zero forcing sets, zero forcing number, k- forcing sets, k-forcing number, forcing number, oriented complete graphs}\\
\indent {\small \textbf{AMS subject classification: 05C69}}

\section{Introduction}
In this paper we discuss the maximum $k$-forcing number over all orientations of a complete graph, which is an interesting case of a more general concept recently introduced by Caro, Davila and Pepper \cite{CaroDavilaPepper}. These concepts generalize both the directed zero forcing number, first introduced in \cite{Hogben2010} and studied in \cite{MinRankDigraphs}, while also expanding recent work on the $k$-forcing number introduced in \cite{Amos2014} and studied further in \cite{Dynamic}. The idea of zero forcing (for simple graphs) was introduced independently in \cite{AIM} and \cite{Burgarth}. In \cite{AIM}, zero forcing was used to bound from below the minimum rank of a graph, or equivalently, to bound from above the maximum nullity of a graph. In \cite{Burgarth}, it is indirectly introduced in relation to a study of control of quantum systems. Additionally, the problem of zero forcing number is closely related to the Power Dominating Set problem, which is motivated by monitoring electric power networks using Kirchoff's Law \cite{power}. Many other papers have been written about zero forcing in recent years (for example \cite{Davila Kenter, LuTang, Yi}). While most of the first papers written were from a linear algebra point of view (\cite{MinRankDigraphs, Edholm, Row}), a fruitful change to a graph theoretic approach, and connection to basic graph parameters such as degree and connectivity, as well as the more general notion of $k$-forcing, was introduced and developed in \cite{Amos2014}, \cite{Dynamic} and \cite{CaroDavilaPepper}. The main point of this paper is to focus the attention on complete graphs and complete $q$-partite graphs, where we already get some interesting results.

Let $G$ be a finite and simple undirected graph with vertex set $V = V(G)$ and edge set $E = E(G)$. We say that $G$ is \emph{oriented} by assigning to each edge $\{u,v\} \in E$ exactly one of the ordered pairs $(u,v)$ and $(v,u)$ -- which we call \emph{arcs}. We call the resulting digraph $D$ an \emph{orientation} of $G$, and say that $D$ is an oriented graph with underlying graph $G$. Let $D$ be an oriented graph with underlying simple graph $G$. If $(u,v)$ is an arc of $D$, then we say that $u$ is \emph{directed towards} $v$, that $v$ is an \emph{out-neighbor} of $u$, and that $u$ is an \emph{in-neighbor} of $v$. Following standard notation: we use $n = n(G), \delta=\delta(G)$ and $\Delta=\Delta(G)$ to denote the order of $G$, the minimum degree of $G$ and the maximum degree of $G$ respectively. A graph with $n=1$ is called a \emph{trivial graph}. If $E = \emptyset$, we say that $G$ is the \emph{empty graph}; otherwise $G$ is a \emph{non-empty graph}. The \emph{degree} of a vertex $v$ is denoted $d(v)$. For any vertex $v$ of $D$, the \emph{out-degree} (resp. \emph{in-degree}) of $v$ is denoted by $d^+(v)$ (resp. $d^-(v)$), and is the number of out-neighbors of $v$ (resp. in-neighbors of $v$). The \emph{minimum out-degree} (resp. in-degree) is denoted $\delta^+ = \delta^+(D)$ (resp. $\delta^- = \delta^-(D)$), and the maximum out-degree (resp. in-degree) is denoted $\Delta^+ = \Delta^+(D)$ (resp. $\Delta^- = \Delta^-(G)$). If every vertex has the same out-degree (resp. in-degree), then $D$ is said to be \emph{out-regular} (resp. \emph{in-regular}). A \emph{directed path} in $D$ is a sequence of vertices $u_1, u_2, \ldots, u_p$ of $D$ such that $(u_i, u_{i+1})$ is an arc of $D$, $1 \leq i \leq p-1$. For terms not defined here, the reader is referred to \cite{West}.

Now we will describe the $k$-forcing process for oriented graphs. Suppose that $D$ is an orientation of $G$, and $S$ is some subset of colored vertices in $D$, all vertices not in $S$ being non-colored. For each positive integer $k$, we define the \emph{$k$-color change rule} as follows: any colored vertex that is directed towards at most $k$ non-colored vertices (has at most $k$ non-colored out-neighbors) forces each of these non-colored vertices to become colored. A colored vertex that forces a non-colored vertex to become colored is said to \emph{$k$-force} that vertex to change color. By the \emph{oriented $k$-forcing process starting from $S \subseteq V$}, we mean the process of first coloring the vertices of $S$, and then iteratively applying the $k$-color change rule as many times as possible. During each step (or iteration) of the oriented $k$-forcing process, all vertices that $k$-force do so simultaneously. If, after termination of the oriented $k$-forcing process, every vertex of $D$ is colored, we say that $S$ is an \emph{oriented $k$-forcing set} (or simply a \emph{$k$-forcing set}) for $D$. The cardinality of a smallest oriented $k$-forcing set for $D$ is called the \emph{oriented $k$-forcing number} of $D$ and is denoted $F_k(D)$. When $k=1$, we will drop the subscript from our notation and write $F(D)$ instead of $F_1(D)$, and this case corresponds to the directed zero forcing number. The maximum oriented $k$-forcing number, over all orientations of $G$, is denoted $\MOF_k(G)$. The minimum oriented $k$-forcing number, over all orientations of $G$, which is denoted $\mof_k(G)$ was also introduced and studied in \cite{CaroDavilaPepper}. These graph invariants turn out to be related to some other well studied graph parameters. For instance, in \cite{CaroDavilaPepper} it is shown that $\MOF_k(G) \geq \alpha(G)$ and, when $k=1$, $\mof_1(G)=\mof(G)=\rho(G)$, where $\alpha(G)$ is the independence number and $\rho(G)$ is the path covering number.   

The remainder of the paper is organized as follows. In Section \ref{MAIN}, we present our results about $\MOF_k(K_n)$. 
In Section \ref{complete q}, we consider $\MOF_k(G)$ when $G$ is a complete $q$-partite graph. In Section \ref{conclusion}, we offer some concluding remarks and acknowledgments.

As a notational convenience, we will use $\log(n)$ (in place of $\log_2(n)$) to denote the base 2 logarithm of $n$, and $k$ will always denote a positive integer.


\section{Main results}\label{MAIN}
In this section, we study the maximum oriented $k$-forcing number for complete graphs. In what follows, we will need to recall that a \emph{transitive orientation} of $D$, with vertices labeled $\{v_1,v_2,\ldots, v_n\}$, is an orientation which satisfies: $v_i$ is directed towards $v_j$ if and only if $i<j$. Also, a \emph{balanced orientation} of $D$ is an orientation satisfying the inequality $|d^+(v) - d^-(v)|\leq 1$, for every vertex $v$. Note that while a complete graph has only one transitive orientation, up to isomorphism, it can have many different non-isomorphic balanced orientations. It turns out that while the forcing number for the transitive orientation is about half the order (seen below), certain kinds of balanced orientations (or nearly balanced) can have forcing numbers that are quite high. 
\begin{prop}\label{transitive}
If $D$ is a transitive orientation of $K_n$, then
\[F_k(D) = \Big\lceil \frac{n}{k+1} \Big\rceil. \]
\end{prop}
\proof Let $D$ be the transitive orientation of $G$, suppose $G$ has $n=q(k+1)+r$ vertices where $0 \leq r < k+1$, $n \geq 2$ and $k<n$. Label the vertices so that $v_i$ has in-degree $n-i$ and out-degree $i-1$. So, for illustration, $v_1$ has in-degree $n-1$ and out-degree $0$, $v_2$ has in-degree $n-2$ and out-degree $1$ and $v_n$ has in-degree $0$and out-degree $n-1$. First we will show that $F_k(D) \leq \lceil\frac{n}{k+1}\rceil$. 

Consider the set $S=\{v_{j(k+1)}\}_{j=1}^q \cup \{v_n\}$, where if $r=0$, $\{v_n\}=\emptyset$. Due to the transitive orientation, no vertex in $S$ with higher label can $k$-force before vertices in $S$ with lower labels. Since $v_{k+1}$ is the lowest labeled vertex in $S$, we start by coloring that vertex. The vertex $v_{k+1}$ has exactly $(k+1)-1=k$ out-neighbors and can $k$-force all of them to change color on the first step of the $k$-forcing process. Once the $k$ out-neighbors of $v_{k+1}$ are colored, the vertex $v_{2(k+1)}$ with out-degree $i-1=2(k+1)-1=2k+1$ can color its $k$ non-colored out-neighbors ($k=1$ of its out-neighbors are already colored). This process continues, with $v_{j(k+1)}$ coloring its remaining non-colored out-neighbors only after $v_{(j-1)(k+1)}$ colors its non-colored out-neighbors until $j=q$. At the last step, if $r=0$ everything is colored, and if $r>0$, the vertex $v_n$ will color the remaining non-colored vertices since there will be at most $(n-1)-q(k+1)=r-1<k+1-1=k$ of them. This shows that $S$ is a $k$-forcing set with $q=\lceil \frac{n}{k+1} \rceil$ vertices if $r=0$ and $q+1=\lceil \frac{q(k+1)+r}{k+1} \rceil = \lceil \frac{n}{k+1} \rceil$ vertices if $r>0$. Thus, in either case we have $F_k(D) \leq \lceil \frac{n}{k+1} \rceil$.

Next we show that $F_k(D) \geq \lceil \frac{n}{k+1} \rceil$. Proceeding by contradiction, assume $F_k(D) < \lceil \frac{n}{k+1} \rceil$ and let $S$ be a smallest $k$-forcing set with $|S| < \lceil \frac{n}{k+1} \rceil$. After coloring each vertex of $S$, and observing that each vertex could $k$-force at most $k$ others, the total number of vertices that end up colored is $|S|+k|S|=|S|(k+1)$. Since $S$ was an oriented $k$-forcing set, every vertex must have been colored so $|S|(k+1) \geq n = q(k+1)+r$. Thus, $|S| \geq q + \frac{r}{k+1}$. However, since $|S| < \lceil \frac{n}{k+1} \rceil = q + \lceil \frac{r}{k+1} \rceil$, and since $|S|$ is an integer, we reach a contradiction and proves the theorem. \qed

The result above leads to the following corollary, which partially supports a conjecture in \cite{CaroDavilaPepper}, namely that $\MOF_k(G) \geq \lceil \frac{n}{k+1} \rceil$, and in particular, $\MOF(G) \geq \lceil \frac{n}{2} \rceil$.
\begin{cor}\label{tran-lower}
For all positive integers $n$, 
\[ \MOF_k(K_n) \geq \Big\lceil \frac{n}{k+1} \Big\rceil. \]
\end{cor}
\proof This follows because $\MOF_k(K_n)$ is at least as much as the oriented $k$-forcing number of the transitive orientation which is $\lceil \frac{n}{k+1} \rceil$, as seen above. \qed

We next recall two results from \cite{CaroDavilaPepper}.
\begin{thm}\label{upper MOF}\cite{CaroDavilaPepper}
Let $G$ be a graph with $n$ vertices and let $D$ be an orientation of $G$ which realizes $\MOF(G)$, so that $F(D)=\MOF(G)$. If $H$ is an induced subgraph of $D$, then 
\[ \MOF(G)\leq F(H)+n-|H|\leq \MOF(H)+n-|H|. \]
\end{thm}

\begin{prop}\label{prop: MOF monotonicity}\cite{CaroDavilaPepper}
If $H$ is any induced subgraph of a graph $G$, then $\MOF_k(G) \geq \MOF_k(H)$.
\end{prop}

Applying these results to complete graphs, we get the following corollary.
\begin{cor}\label{one_more}
If $n$ is a positive integer, then $\MOF(K_n) \leq \MOF(K_{n+1}) \leq \MOF(K_n) +1$.
\end{cor}
\proof The lower bound comes from Proposition \ref{prop: MOF monotonicity}. For the upper bound, let $H$ be an induced $K_n$ inside of a $K_{n+1}$. From Theorem \ref{upper MOF},
\[ \MOF(K_{n+1}) \leq \MOF(K_n) + (n+1) - |H| = \MOF(K_n) + (n+1) - n = \MOF(K_n)+1. \]
\qed

\begin{table}[tbp]
    \centering
    \begin{tabular}{|c||c|c|c|c|c|c|c|c|c|c|c|c|c|c|c|c|c|c|c|}
         \hline
         Order &  3 & 4 & 5 & 6 &7&8&9&10&11&12&13&14&15&16&17&18&19 \\ \hline
         MOF & 2&2&3&3&4&5&6&6&7&8&8&9&10&10&11&12&13 \\ \hline
    \end{tabular}
\caption{Values of $MOF(K_n)$ found from a computer program employing various theorems from \cite{CaroDavilaPepper}, as well as Corollary \ref{one_more}, and with many processors running in parallel. }
    \label{tab:my_label}
\end{table}

Since $\MOF(K_n)$ can grow by at most one as $n$ grows by one, Corollary~\ref{one_more} considerably speeds up any attempt to find exactly the values of $\MOF(K_n)$. Further speed ups to any computation of $\MOF(K_n)$ come from other theorems in \cite{CaroDavilaPepper}, in particular the Reversal Theorem. Namely, that the forcing number of an orientation of a graph is equal to the forcing number of its reversal. These ideas, and others, were used to write a computer program to find the exact value of $\MOF(K_n)$ for $n < 20$. The results of this can be seen in Table~\ref{tab:my_label}. We are now ready to present our main results, which are lower bounds for $\MOF(K_n)$.

\begin{thm}\label{3-4}
If $G$ is a graph with order $n$, then
\[
\MOF(K_n) \geq \frac{3n-9}{4}. 
\]
\end{thm}
\proof As can be seen from inspection of Table 1 (below), this theorem is true for all values of $n \leq 10$ (achieving equality when $n=6$ and $n=10$ if we consider $\lceil \frac{3n-9}{4} \rceil$ since $\MOF(K_n)$ is an integer). Let us then assume, without loss of generality, that $n \geq 10$. 

Let $q$ be the largest odd integer such that $n=2q+r$ where $0 \leq r < q$. It can be readily seen then that $r \in \{0,1,2,3\}$, since if $r \geq 4$, $q$ was not the largest odd integer satisfying the equation. Hence, it is implied that $0 \leq r \leq 3$. 

We consider the following orientation $D$ of $K_n$. Partition the vertices into $q-r$ sets of order $2$ and $r$ sets of order $3$. Label the $q-r$ sets of order $2$ as $\{A_1,A_2,\ldots,A_{q-r}\}$ and label the $r$ sets of order $3$ as $\{A_{q-r+1},A_{q-r+2},\ldots,A_q\}$. Now split these sets into two nearly equal halves, with one having $\frac{q-1}{2}$ parts and the other having $\frac{q+1}{2}$ parts. Let the lowered labeled sets be in the smaller of these groups, so that each of the sets $\{A_1,A_2,\ldots,A_{\frac{q-1}{2}}\}$ has order 2. Since $n \geq 10$, we know that $q \geq 5$, and this implies that all $r$ extra vertices are in the higher labeled group. To ease the notation, let $X$ denote the lowest $\frac{q-1}{2}$ labeled sets and let $Y$ denote the highest $\frac{q+1}{2}$ labeled sets. Now, give each of these sets $A_i$ the transitive orientation with respect to the other vertices in that set. Consider the sets themselves as vertices in larger odd order graph, and give that graph a balanced orientation in the following way. Each vertex in $A_i$ is joined to each vertex in each of the next $\frac{q-1}{2}$ highest labeled sets (wrapping around again when we get past $A_q$). So for example, if $q=7$, then each vertex of $A_2$ is joined to each vertex of $A_3$, $A_4$, and $A_5$. 

Now we are ready to consider how many vertices need to be colored to have a chance at forcing the whole graph. In order for any vertex in $Y$ to force any other vertex to change color, the initial set of colored vertices must be at least as large as $|X|-1$. Considering the set $Y$ as a separate oriented complete graph, we discover that it has the transitive orientation. Thus, according to Theorem \ref{transitive}, in order for that set to be colored, once $X$ is colored, we need at least $\lceil \frac{|Y|}{2} \rceil$ many vertices to be initially colored. Taken together,
\[
\MOF(K_n) \geq F(D) \geq |X|-1 + \Big\lceil \frac{|Y|}{2} \Big\rceil \geq 2\Big(\frac{q-1}{2}\Big)-1 + \frac{2(\frac{q+1}{2})+r}{2} = \frac{3n-6-r}{4}.
\]
Finally, since $r \leq 3$, the result follows and the theorem is proven. \qed


\begin{figure}[htb]\label{fig1}
\begin{center}
\begin{tikzpicture}[scale=.8,style=thin,x=1cm,y=1cm, =>stealth]
\def\vr{2.5pt} 

\path (0,3.25) coordinate (x1);
\path (-1,3) coordinate (x2);
\path (-2,2) coordinate (x3);
\path (-2,.5) coordinate (x4);
\path (-.75,-.5) coordinate (x5);
\path (2,.5) coordinate (x7);
\path (.75,-.5) coordinate (x6);
\path (2, 2) coordinate (x8);
\path (1, 3) coordinate (x9);

\path (0.75,1.3) coordinate (w);
%
\draw[black, arrows={->[line width=.25pt,black,length=3.5mm,width=1.25mm]}] (x1) -- (x2);
\draw[black, arrows={->[line width=.25pt,black,length=3.5mm,width=1.25mm]}]  (x1) -- (x3);
\draw[black, arrows={->[line width=.25pt,black,length=3.5mm,width=1.25mm]}]  (x1) -- (x4);
\draw[black, arrows={->[line width=.25pt,black,length=3.5mm,width=1.25mm]}] (x1) -- (x5);
\draw[black, arrows={->[line width=.25pt,black,length=3.5mm,width=1.25mm]}] (x1) -- (x6);
\draw[black, arrows={->[line width=.25pt,black,length=3.5mm,width=1.25mm]}] (x1) -- (x7);
\draw[black, arrows={->[line width=.25pt,black,length=3.5mm,width=1.25mm]}] (x1) -- (x8);
\draw[black, arrows={->[line width=.25pt,black,length=3.5mm,width=1.25mm]}]  (x1) -- (x9);

\draw[black, arrows={->[line width=.25pt,black,length=3.5mm,width=1.25mm]}] (x2) -- (x3);
\draw[black, arrows={->[line width=.25pt,black,length=3.5mm,width=1.25mm]}] (x2) -- (x4);
\draw[black, arrows={->[line width=.25pt,black,length=3.5mm,width=1.25mm]}] (x2) -- (x5);
\draw[black, arrows={->[line width=.25pt,black,length=3.5mm,width=1.25mm]}] (x2) -- (x6);
\draw[black, arrows={->[line width=.25pt,black,length=3.5mm,width=1.25mm]}] (x2) -- (x7);
\draw[black, arrows={->[line width=.25pt,black,length=3.5mm,width=1.25mm]}] (x8) -- (x2);
\draw[black, arrows={->[line width=.25pt,black,length=3.5mm,width=1.25mm]}] (x9) -- (x2);

\draw[black, arrows={->[line width=.25pt,black,length=3.5mm,width=1.25mm]}] (x3) -- (x4);
\draw[black, arrows={->[line width=.25pt,black,length=3.5mm,width=1.25mm]}] (x3) -- (x5);
\draw[black, arrows={->[line width=.25pt,black,length=3.5mm,width=1.25mm]}] (x3) -- (x6);
\draw[black, arrows={->[line width=.25pt,black,length=3.5mm,width=1.25mm]}] (x3) -- (x7);
\draw[black, arrows={->[line width=.25pt,black,length=3.5mm,width=1.25mm]}] (x8) -- (x3);
\draw[black, arrows={->[line width=.25pt,black,length=3.5mm,width=1.25mm]}] (x9) -- (x3);

\draw[black, arrows={->[line width=.25pt,black,length=3.5mm,width=1.25mm]}] (x4) -- (x5);
\draw[black, arrows={->[line width=.25pt,black,length=3.5mm,width=1.25mm]}] (x4) -- (x6);
\draw[black, arrows={->[line width=.25pt,black,length=3.5mm,width=1.25mm]}] (x4) -- (x7);
\draw[black, arrows={->[line width=.25pt,black,length=3.5mm,width=1.25mm]}] (x4) -- (x8);
\draw[black, arrows={->[line width=.25pt,black,length=3.5mm,width=1.25mm]}] (x4) -- (x9);

\draw[black, arrows={->[line width=.25pt,black,length=3.5mm,width=1.25mm]}] (x5) -- (x6);
\draw[black, arrows={->[line width=.25pt,black,length=3.5mm,width=1.25mm]}] (x5) -- (x7);
\draw[black, arrows={->[line width=.25pt,black,length=3.5mm,width=1.25mm]}] (x5) -- (x8);
\draw[black, arrows={->[line width=.25pt,black,length=3.5mm,width=1.25mm]}] (x5) -- (x9);

\draw[black, arrows={->[line width=.25pt,black,length=3.5mm,width=1.25mm]}] (x6) -- (x7);
\draw[black, arrows={->[line width=.25pt,black,length=3.5mm,width=1.25mm]}] (x6) -- (x8);
\draw[black, arrows={->[line width=.25pt,black,length=3.5mm,width=1.25mm]}] (x6) -- (x9);

\draw[black, arrows={->[line width=.25pt,black,length=3.5mm,width=1.25mm]}] (x8) -- (x7);
\draw[black, arrows={->[line width=.25pt,black,length=3.5mm,width=1.25mm]}] (x9) -- (x7);

\draw[black, arrows={->[line width=.25pt,black,length=3.5mm,width=1.25mm]}] (x8) -- (x9);
\draw (x1) [fill=black]circle (\vr);
\draw (x2) [fill=black]circle (\vr);
\draw (x3) [fill=black] circle (\vr);
\draw (x4) [fill=black] circle (\vr);
\draw (x5) [fill=black] circle (\vr);
\draw (x6) [fill=white] circle (\vr);
\draw (x7) [fill=white] circle (\vr);
\draw (x8) [fill=white] circle (\vr);
\draw (x9) [fill=black] circle (\vr);

\draw (0,-1.5) node {The graph $K_9$};

\draw[anchor = south] (x1) node {{\small $x_1$}};
\draw[anchor = east] (x2) node {{\small $x_2$}};
\draw[anchor = east] (x3) node {{\small $x_3$}};
\draw[anchor = east] (x4) node {{\small $x_4$}};
\draw[anchor = north] (x5) node {{\small $x_5$}};
\draw[anchor = north] (x6) node {{\small $x_6$}};
\draw[anchor = west] (x7) node {{\small $x_7$}};
\draw[anchor = west] (x8) node {{\small $x_8$}};
\draw[anchor = west] (x9) node {{\small $x_9$}};

\path (8,6.75) coordinate (x1);
\path (6.25,6.25) coordinate (x2);
\path (4.75,5) coordinate (x3);
\path (4,3.25) coordinate (x4);
\path (4.75,1.5) coordinate (x5);
\path (9.5,0) coordinate (x7);
\path (6.25,0) coordinate (x6);
\path (11.25, 1.5) coordinate (x8);
\path (12., 3.25) coordinate (x9);
\path (11.25, 5) coordinate (x10);
\path (9.75, 6.25) coordinate (x11);

\draw[black, arrows={->[line width=.25pt,black,length=3.5mm,width=1.25mm]}] (x1) -- (x2);
\draw[black, arrows={->[line width=.25pt,black,length=3.5mm,width=1.25mm]}] (x1) -- (x3);
\draw[black, arrows={->[line width=.25pt,black,length=3.5mm,width=1.25mm]}] (x1) -- (x4);
\draw[black, arrows={->[line width=.25pt,black,length=3.5mm,width=1.25mm]}] (x1) -- (x5);
\draw[black, arrows={->[line width=.25pt,black,length=3.5mm,width=1.25mm]}] (x1) -- (x6);
\draw[black, arrows={->[line width=.25pt,black,length=3.5mm,width=1.25mm]}] (x1) -- (x7);
\draw[black, arrows={->[line width=.25pt,black,length=3.5mm,width=1.25mm]}] (x1) -- (x8);
\draw[black, arrows={->[line width=.25pt,black,length=3.5mm,width=1.25mm]}] (x1) -- (x9);
\draw[black, arrows={->[line width=.25pt,black,length=3.5mm,width=1.25mm]}] (x1) -- (x10);
\draw[black, arrows={->[line width=.25pt,black,length=3.5mm,width=1.25mm]}] (x1) -- (x11);

\draw[black, arrows={->[line width=.25pt,black,length=3.5mm,width=1.25mm]}] (x2) -- (x3);
\draw[black, arrows={->[line width=.25pt,black,length=3.5mm,width=1.25mm]}] (x2) -- (x4);
\draw[black, arrows={->[line width=.25pt,black,length=3.5mm,width=1.25mm]}] (x2) -- (x5);
\draw[black, arrows={->[line width=.25pt,black,length=3.5mm,width=1.25mm]}] (x2) -- (x6);
\draw[black, arrows={->[line width=.25pt,black,length=3.5mm,width=1.25mm]}] (x2) -- (x7);
\draw[black, arrows={->[line width=.25pt,black,length=3.5mm,width=1.25mm]}] (x2) -- (x8);
\draw[black, arrows={->[line width=.25pt,black,length=3.5mm,width=1.25mm]}] (x2) -- (x9);
\draw[black, arrows={->[line width=.25pt,black,length=3.5mm,width=1.25mm]}] (x2) -- (x10);
\draw[black, arrows={->[line width=.25pt,black,length=3.5mm,width=1.25mm]}] (x2) -- (x11);

\draw[black, arrows={->[line width=.25pt,black,length=3.5mm,width=1.25mm]}] (x3) -- (x4);
\draw[black, arrows={->[line width=.25pt,black,length=3.5mm,width=1.25mm]}] (x3) -- (x5);
\draw[black, arrows={->[line width=.25pt,black,length=3.5mm,width=1.25mm]}] (x3) -- (x6);
\draw[black, arrows={->[line width=.25pt,black,length=3.5mm,width=1.25mm]}] (x3) -- (x7);
\draw[black, arrows={->[line width=.25pt,black,length=3.5mm,width=1.25mm]}] (x3) -- (x8);
\draw[black, arrows={->[line width=.25pt,black,length=3.5mm,width=1.25mm]}] (x3) -- (x9);
\draw[black, arrows={->[line width=.25pt,black,length=3.5mm,width=1.25mm]}] (x3) -- (x10);
\draw[black, arrows={->[line width=.25pt,black,length=3.5mm,width=1.25mm]}] (x11) -- (x3);

\draw[black, arrows={->[line width=.25pt,black,length=3.5mm,width=1.25mm]}] (x4) -- (x5);
\draw[black, arrows={->[line width=.25pt,black,length=3.5mm,width=1.25mm]}] (x4) -- (x6);
\draw[black, arrows={->[line width=.25pt,black,length=3.5mm,width=1.25mm]}] (x4) -- (x7);
\draw[black, arrows={->[line width=.25pt,black,length=3.5mm,width=1.25mm]}] (x4) -- (x8);
\draw[black, arrows={->[line width=.25pt,black,length=3.5mm,width=1.25mm]}] (x4) -- (x9);
\draw[black, arrows={->[line width=.25pt,black,length=3.5mm,width=1.25mm]}] (x10) -- (x4);
\draw[black, arrows={->[line width=.25pt,black,length=3.5mm,width=1.25mm]}] (x11) -- (x4);

\draw[black, arrows={->[line width=.25pt,black,length=3.5mm,width=1.25mm]}] (x5) -- (x6);
\draw[black, arrows={->[line width=.25pt,black,length=3.5mm,width=1.25mm]}] (x5) -- (x7);
\draw[black, arrows={->[line width=.25pt,black,length=3.5mm,width=1.25mm]}] (x5) -- (x8);
\draw[black, arrows={->[line width=.25pt,black,length=3.5mm,width=1.25mm]}] (x9) -- (x5);
\draw[black, arrows={->[line width=.25pt,black,length=3.5mm,width=1.25mm]}] (x10) -- (x5);
\draw[black, arrows={->[line width=.25pt,black,length=3.5mm,width=1.25mm]}] (x11) -- (x5);

\draw[black, arrows={->[line width=.25pt,black,length=3.5mm,width=1.25mm]}] (x6) -- (x7);
\draw[black, arrows={->[line width=.25pt,black,length=3.5mm,width=1.25mm]}] (x6) -- (x8);
\draw[black, arrows={->[line width=.25pt,black,length=3.5mm,width=1.25mm]}] (x6) -- (x9);
\draw[black, arrows={->[line width=.25pt,black,length=3.5mm,width=1.25mm]}] (x6) -- (x10);
\draw[black, arrows={->[line width=.25pt,black,length=3.5mm,width=1.25mm]}] (x6) -- (x11);

\draw[black, arrows={->[line width=.25pt,black,length=3.5mm,width=1.25mm]}] (x7) -- (x8);
\draw[black, arrows={->[line width=.25pt,black,length=3.5mm,width=1.25mm]}] (x7) -- (x9);
\draw[black, arrows={->[line width=.25pt,black,length=3.5mm,width=1.25mm]}] (x7) -- (x10);
\draw[black, arrows={->[line width=.25pt,black,length=3.5mm,width=1.25mm]}] (x7) -- (x11);

\draw[black, arrows={->[line width=.25pt,black,length=3.5mm,width=1.25mm]}] (x8) -- (x9);
\draw[black, arrows={->[line width=.25pt,black,length=3.5mm,width=1.25mm]}] (x8) -- (x10);
\draw[black, arrows={->[line width=.25pt,black,length=3.5mm,width=1.25mm]}] (x8) -- (x11);

\draw[black, arrows={->[line width=.25pt,black,length=3.5mm,width=1.25mm]}] (x10) -- (x9);
\draw[black, arrows={->[line width=.25pt,black,length=3.5mm,width=1.25mm]}] (x11) -- (x9);

\draw[black, arrows={->[line width=.25pt,black,length=3.5mm,width=1.25mm]}] (x10) -- (x11);

\draw (8,-1.5) node {The graph $K_{11}$};
\draw (x1) [fill=black]circle (\vr);
\draw (x2) [fill=black]circle (\vr);
\draw (x3) [fill=white] circle (\vr);
\draw (x4) [fill=black] circle (\vr);
\draw (x5) [fill=black] circle (\vr);
\draw (x6) [fill=black] circle (\vr);
\draw (x7) [fill=black] circle (\vr);
\draw (x8) [fill=white] circle (\vr);
\draw (x9) [fill=white] circle (\vr);
\draw (x10) [fill=white] circle (\vr);
\draw (x11) [fill=black] circle (\vr);

\draw[anchor = south] (x1) node {{\small $x_1$}};
\draw[anchor = east] (x2) node {{\small $x_2$}};
\draw[anchor = east] (x3) node {{\small $x_3$}};
\draw[anchor = east] (x4) node {{\small $x_4$}};
\draw[anchor = north] (x5) node {{\small $x_5$}};
\draw[anchor = north] (x6) node {{\small $x_6$}};
\draw[anchor = north] (x7) node {{\small $x_7$}};
\draw[anchor = west] (x8) node {{\small $x_8$}};
\draw[anchor = west] (x9) node {{\small $x_9$}};
\draw[anchor = west] (x10) node {{\small $x_{10}$}};
\draw[anchor = west] (x11) node {{\small $x_{11}$}};

\end{tikzpicture}
\end{center}
\vskip -0.5 cm
\caption{Orientations of $K_9$ and $K_{11}$, with minimum forcing sets highlighted.} \label{f:fig1}
\end{figure}
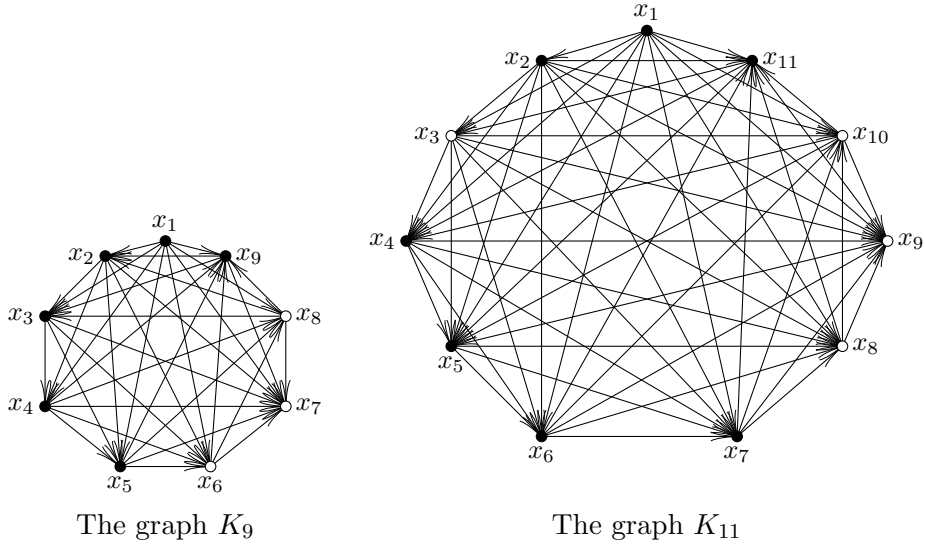


\begin{lem}\label{tech_lem}
If $n =3p+r \geq 9$, where $r \in \{0,1,2\}$, then $\MOF(K_n) \geq p + 2\MOF(K_p) - 3$.
\end{lem}

\proof Assume $n=3p+r \geq 9$ with $r \in \{0,1,2\}$. Partition $K_n$ into three parts whose orders are as close to equal as possible. Label these three parts as $V_1$, $V_2$ and $V_3$. Let $D$ denote the orientation of $K_p$ which realizes $\MOF(K_p)$, and let $D_i$ denote the orientation of $V_i$ which realizes $MOF(K_{|V_i|})$ for $i \in \{1,2,3\}$. Now, orient all other edges of the graph as follows. Each vertex in $V_1$ is directed toward each vertex of $V_2$, each vertex of $V_2$ is directed towards each vertex of $V_3$, and each vertex of $V_3$ is directed towards each vertex of $V_1$. Call the completed orientation of $K_n$ thus created $D^*$, and let $S$ be a minimum oriented forcing set of $D^*$. Finally, for each $i \in \{1,2,3\}$, let $S_i = S \cap V_i$. 

Not all vertices are originally colored, so there must be a vertex which forces at the first step of the forcing process. Let $v$ be such a vertex and assume that $v \in V_1$. There are two main cases to consider, that $v$ forces a vertex in $V_2$ or that $v$ forces a vertex in $V_1$. 

First, suppose $v$ and forces a vertex in $V_2$. This is only possible if all but one of the vertices of $V_2$ are already colored, implying $|S_2|=|V_2|-1$. Next, in order for $V_1$ itself to be fully colored, either $S_1$ itself is an oriented forcing set of $V_1$, which implies $|S_1| \geq F(D_1)$, or the last non-colored vertex from $V_1$ is forced by a vertex in $V_3$. This later situation is only possible if all but one vertex in $V_1$ is already colored, which implies $S_1$ was able to color all but one vertex of $V_1$ so that $|S_1| \geq F(D_1)-1$. Finally, in order now for $V_3$ to be fully colored, either $S_3$ itself was a forcing set of $V_3$, which implies $|S_3| \geq F(D_3)$, or the last non-colored vertex of $V_3$ is forced by a vertex in $V_2$. This later situation is only possible if all but one vertex of $V_3$ is already colored, which implies $S_3$ was able to color all but one vertex in $V_3$ so that $|S_3| \geq F(D_3) -1$. Hence, summing the parts, we get:

\[
|S|=|S_1|+|S_2|+|S_3| \geq (F(D_1)-1)+(|V_2|-1)+(F(D_3)-1) = |V_2| + F(D_1)+F(D_3)-3.
\]

Second, suppose $v$ and forces a vertex in $V_1$. This is clearly only possible if all of $V_2$ is already colored, which implies $|S_2|=|V_2|$. Now, the argument repeats as in the preceding paragraph. In order for $V_1$ to be fully colored, $|S_1| \geq F(D_1)-1$ and in order for $V_3$ to be fully colored, $|S_3| \geq F(D_3)-1$. Hence, summing the parts, we get:

\[
|S|=|S_1|+|S_2|+|S_3| \geq (F(D_1)-1)+(|V_2|)+(F(D_3)-1) > |V_2| + F(D_1)+F(D_3)-3.
\]
Thus, in either case we have, 

\begin{equation}\label{A}
F(D^*)=|S| \geq |V_2| + F(D_1)+F(D_3)-3.
\end{equation}

If we assume $v \in V_2$ instead of $V_1$, the argument above could be repeated and we would arrive at the inequality,

\begin{equation}\label{B}
F(D^*)=|S| \geq |V_3| + F(D_1)+F(D_2)-3.
\end{equation}

If we assume $v \in V_3$ instead of $V_1$, the argument above could be repeated and we would arrive at the inequality,

\begin{equation}\label{C}
F(D^*)=|S| \geq |V_1| + F(D_2)+F(D_3)-3.
\end{equation}

To conclude, we make use of the facts that for each $i \in \{1,2,3\}$, $|V_i| \geq p$ and $F(D_i) \geq F(D)=\MOF(K_p)$, to bound from below each of Inequalities \ref{A}, \ref{B} and \ref{C} by $p + 2MOF(K_p)-3$. Therefore, $\MOF(K_n) \geq F(D^*) \geq p + 2\MOF(K_p)-3$ as claimed. \qed

Theorem \ref{3-4} and Lemma \ref{tech_lem} can be used together to help us prove the following main result, namely that $\MOF(K_n)$ is asymptotically equal to $n$.

\begin{thm}\label{log} For all positive integers $n \geq 2$, 
\[
\MOF(K_n) \geq n - \frac{2n}{\log(n)}.
\]
\end{thm}
\proof Proceeding by mathematical induction, notice that for all values of $n$ in the range, $2 \leq n \leq 202$, the results follows from Theorem \ref{3-4} since,
\[
\MOF(K_n) \geq \frac{3n-9}{4} > n - \frac{2n}{\log(n)},
\]
when $2 \leq n \leq 202$ as seen from calculation and inspection. This settles our base case. Assume now that the theorem is true for all complete graphs with smaller orders than $K_n$, with $n \geq 203$. We will show this implies it is also true for $K_n$. Let $n=3p+r$, with $r \in \{1,2,3\}$. Now, from Lemma \ref{tech_lem} we know that,
\[
\MOF(K_n) \geq p + 2\MOF(K_p) - 3.
\]
Since $p < n$, from our inductive assumption we know that $\MOF(K_p) \geq p - \frac{2p}{\log(p)}$. Together, this implies,
\[
\MOF(K_n) \geq 3p - \frac{4p}{\log(p)} -3.
\]
It remains to show that $3p - \frac{4p}{\log(p)} -3 \geq n - \frac{2n}{\log(n)}$. Replacing $p$ by $\frac{n-r}{3}$, and rearranging the terms, this is equivalent to showing that,
\[
\frac{2n}{\log(n)} \geq \frac{\frac{4}{3}(n-r)}{\log(\frac{n-r}{3})} + r + 3.
\]
Finally, since this last inequality is true for all $n \geq 117$, since the function,
\[f(n)= \frac{2n}{\log(n)} - \frac{\frac{4}{3}(n-r)}{\log(\frac{n-r}{3})} - r - 3,
\]
is never negative for $r$ and $n$ in the ranges given, as can be seen using standard techniques from calculus and algebra, and we already assumed $n\geq 202$. The general result now follows by induction and the theorem is proven. \qed

From this we easily deduce that $\frac{\MOF(K_n)}{n} \rightarrow 1$ as $n \rightarrow \infty$. Furthermore, when combined with Proposition \ref{prop: MOF monotonicity}, we get the following corollary. Recall that the \emph{clique number} of $G$, written $\omega(G)$, is the cardinality of a largest induced complete graph in $G$.

\begin{cor}\label{omega}
If $G$ is a graph with order n, then
\[ \MOF(G) \geq \omega(G) - \frac{2\omega(G)}{\log(\omega(G))}. \]
\end{cor}
\proof Let $H$ be a largest complete subgraph of $G$. From Proposition \ref{prop: MOF monotonicity}, $\MOF(G) \geq \MOF(H)$. Now, since $H$ is a complete graph of order $\omega(G)$, the inequality follows from Theorem \ref{log}. \qed

To conclude this section, we recall one more result from \cite{CaroDavilaPepper}.

\begin{cor}\cite{CaroDavilaPepper}
If $G$ is a graph with order $n$, then 
\[ \MOF(G) \leq n - \frac{\log(\omega(G))}{2}. \]
\end{cor}

Thus, taken together with the observation that $\omega(K_n) = n$, we find the following.

\begin{cor}\label{omega2}
For all positive integers $n$, 
\[ n - \frac{2n}{\log(n)} \leq \MOF(K_n) \leq n - \frac{\log(n)}{2}. \]
\end{cor}

\section{The maximum oriented $k$-forcing number for complete $q$-partite graphs} \label{complete q}


In this section, we extend our investigation to from complete graphs to complete $q$-partite graphs. Recall that a graph $G$ is \emph{$q$-partite} if its vertex set can be partitioned into $q \geq 2$ independent sets. The independent sets are called \emph{parts} and, if $G$ is $q$-partite, every edge in $G$ has its two incident vertices in different parts. We say that $G$ is a \emph{complete $q$-partite graph} if $G$ is $q$-partite with every possible edge between vertices in different parts.
\begin{thm}\label{q-part}
If $G$ is a complete $q$-partite graph and $n_1 \geq n_2 \geq \ldots \geq n_q$ denote the orders of the partite sets, then
\[
\MOF_k(G) \geq n_1 + \sum_{i=2}^q \max{\{n_i-k,0\}}.
\]
\end{thm}
\proof Label the partite sets as $A_1, A_2, \ldots, A_q$, labeled so that larger parts have smaller labels. That is, if $i < j$, then $|A_i| \geq |A_j|$. Moreover, let $n_i=|A_i|$. Create the orientation $D$ of the edges of $G$ by directing vertices from parts with smaller labels towards parts with larger labels. That is, if $u \in A_i$ and $v \in A_j$, with $i < j$, then $(u,v)$ is an arc in $D$ (this is called the transitive orientation). Now, in $D$, all vertices from $A_1$ must be in any oriented $k$-forcing set. Moreover, with $i < j$, vertices from $A_i$ could only $k$-force vertices from $A_j$ if $\max{\{n_i-k,0\}}$ vertices from $A_j$ were already colored. Since this is true for all pairs $i$ and $j$, we get the following lower bound on the oriented $k$-forcing number of $D$,
\[
F_k(D) \geq n_1 + \sum_{i=2}^q \max{\{n_i-k,0\}}.
\]
Finally, since $\MOF_k(G) \geq F_k(D)$, the proof is complete. \qed

\begin{cor}
If $G$ is a complete $q$-partite graph, with $q \geq 2$, and $k$ is a positive integer, then
\[
\MOF_k(G) \geq n - k(q - 1).
\]
\end{cor}
\proof
Let $G$ be a complete $q$-partite graph with parts $A_1, A_2, \ldots, A_q$. Set $|A_i|=n_i$ and without loss of generality, assume $n_1 \geq n_2 \geq \ldots \geq n_q$. From Theorem \ref{q-part} above, together with the fact that $\max{\{n_i-k,0\}} \geq n_i-k$, we get the following chain of inequalities;
\[
\MOF_k(G) \geq n_1 + \sum_{i=2}^q \max{\{n_i-k,0\}} \geq n_1 + \sum_{i=2}^q (n_i-k) = n - k(q-1),
\]
which completes the proof. \qed

Specifying that $k=1$, we arrive at the result below.
\begin{cor}\label{partite lower}
If $G$ is a complete $q$-partite graph, with $q \geq 2$, then
\[
\MOF(G) \geq n - q + 1.
\]
This inequality is sharp when $q=2$ or when $q=3$ and each part has at least two vertices.
\end{cor}
\proof The inequality comes from substituting $k=1$ into the above corollary. To see that equality holds when $q=2$, the complete bipartite case, we first note that $\MOF(G) \leq n-1$ for any non-empty graph $G$. This is because there is always a vertex $v$ with in-degree at least one in such cases, and the set of all vertices other than $v$ is a forcing set. To see that $\MOF(G) \geq n-1$, let $A$ and $B$ be the two parts and direct all edges from $A$ to $B$. Now each vertex of $A$ is necessarily in any forcing set and nothing can be forced unless at least $|B|-1$ vertices from $B$ are included. Therefore, $\MOF(G)=n-1=n-q+1$, as claimed.

Next we show that the inequality is sharp when $q=3$ and each part has at least two vertices. Let $A,B,$ and $C$ be the three parts with cardinalities $a \geq b \geq c \geq 2$ respectively. It is sufficient to show that $\MOF(G) \leq n-q+1=n-2 =a+b+c-2$, since the same lower bound is already established (using transitive orientation). To this end, let $D$ be any orientation of $G$. We first show that there must be two vertices, $u$ and $v$, in two different parts, with in-degree at least one. To see this, suppose there was a vertex $w$ such that $d^+(w)=d(w)$. In this case, each vertex in the two parts not containing $w$ have in-degree at least one and we are done. On the other hand, if no such $w$ exists, then all vertices have $d^+(w)<d(w)$ and consequently have in-degree at least one, and we are done. So, let $u$ and $v$ be vertices in different parts such that $d^-(u)\geq 1$ and $d^-(v)\geq 1$. Without loss of generality, assume $u\in B$ and $v\in C$. Let $u^* \in B$ and $v^* \in C$ be vertices from $B$ and $C$ respectively, different from $u$ and $v$ (since each part has at least two vertices). Suppose, without loss of generality, $(v,u)$ is an arc. This implies there is a vertex $w$ such that $(w,v)$ is an arc since $d^-(v)\geq 1$. If both $(v,u^*)$ and $(u,v^*)$ are arcs, the set $V \setminus \{v^*,u^*\}$ is a forcing set of order $n-2$ since on the first step of the forcing process, $u$ forces $v^*$ and $v$ forces $u^*$. Otherwise, either $(v^*,u)$ or $(u^*,v)$ is an arc. Suppose $(v^*,u)$ is an arc, and consider the set $V \setminus \{v,u\}$. This is a forcing set of order $n-2$, since $v^*$ forces $u$ and then $w$ forces $v$. Suppose $(u^*,v)$ is an arc, and consider the set $V \setminus \{v,u\}$. This is a forcing set of order $n-2$, since $u^*$ forces $v$ and then $v$ forces $u$.  Therefore, $F(D)\leq n-2$ for any orientation $D$. Consequently, $\MOF(G) \leq n-q+1=n-2 =a+b+c-2$, completing the proof.
\qed

For $q \geq 9$, and each of the $q$ parts has at least three vertices, we have examples where $\MOF(G) > n-q+1$, so the inequality is not sharp for large values of $q$. We do not know the situation for $4 \leq q \leq 8$.

If $G$ is a complete $q$-partite graph, then $G$ contains a subgraph isomorphic to $K_q$. This observation together with Theorem~\ref{upper MOF} shows that the maximum oriented forcing number of a $q$-partite graph can be bounded from above in terms of $\MOF(K_q)$. 
In particular, together with Corollary~\ref{partite lower}, we have the following. 
\begin{cor}
If $G$ is a complete $q$-partite graph with order $n$, and $q\ge 2$, then
\begin{equation*}
n - q + 1 \le \MOF(G) \le n - q + \MOF(K_q). 
\end{equation*}
\end{cor}

\section{Concluding remarks and acknowledgments}\label{conclusion}
In this paper we have given a detailed study of the maximum $k$-forcing number over all orientations of complete graphs and complete $q$-partite graphs. However, our focus was primarily on the interesting case when $k = 1$. We highlight that in Corollary~\ref{omega2}, we have the lower bound $\MOF(K_n) \geq n-\frac{2n}{\log(2)}$, and the upper bound $\MOF(K_n) \leq n-\frac{\log(n)}{2}$. It remains to be seen which of these bounds is closer to the truth, and we pose this formally with the following problem.
\begin{prob}
Which of the bounds presented in Corollary~\ref{omega2} is a better approximation to $\MOF(K_n)$?
\end{prob}

Finally, we would like to thank several individuals -- and former students of the second author -- for their help in working on this paper: David Amos for many fruitful conversations about early results; Randy Davila for his great help in preparing and reviewing the manuscript; and Mobeen Azhar for his work in writing and developing the computer program which generated the results from Table~\ref{tab:my_label}.  Their contributions and enthusiasm helped provide the necessary motivation to make this paper possible.

\end{document}